\pdfoutput=1
\pdfoutput=1
\pdfoutput=1
\pdfoutput=1
\pdfoutput=1
\documentclass[a4paper,11 pt]{article}
\usepackage{calc}
\usepackage[all]{xy}
\usepackage[centertags]{amsmath}
\usepackage{latexsym}
\usepackage{amsfonts}
\usepackage{graphicx}
\usepackage{tikz}
\usepackage{cases}
\usepackage{amssymb}
\usepackage{amsthm}
\usepackage{color}
\usepackage{fancyhdr}
\usepackage[dvips]{epsfig}
\usepackage{newlfont}
\usepackage[latin1]{inputenc}
\usepackage[latin1]{inputenc}
\usepackage[english,french]{babel}
\usepackage{graphicx}
\usepackage{t1enc}
\usepackage[english,french]{babel}
\usepackage{fancybox}
\usepackage{graphicx}
\usepackage{t1enc}
\usepackage[french2]{minitoc}
\usepackage{mathrsfs}
\usepackage{amsfonts}
\usepackage{wasysym}
\usepackage{hyperref}
\allowdisplaybreaks
\usepackage{float}
\usepackage{pdfpages}
\usepackage{geometry}

\fancypagestyle{plain}{ \fancyhead{}

	\rhead{\textbf{}} \rfoot{\footnotesize{\textsf{\tiny }}}
	\cfoot{\footnotesize{\textbf{ \thepage}}}} \hfuzz2pt
\newlength{\defbaselineskip}
\numberwithin{equation}{section} 
\newtheorem{theorem}{Theorem}[section]

\newtheorem{lemma}[theorem]{Lemma}
\newtheorem{example}{Example}[section]
\newtheorem{proposition}[theorem]{Proposition}

\newtheorem{remark}{Remark}[section]

\newcommand{\Z}{{\mathbb Z}}

 \makeatletter
\newcommand{\thechapterwords}
{ \ifcase \thechapter\or 1\or 2\or 3\or 4\or 5\or
	6\or 7\or 8\or 9\or 10\or 11\fi}
\def\thickhrulefill{\leavevmode \leaders \hrule height 2ex \hfill \kern \z@}
\def\@makechapterhead#1{%
	\vspace*{15\p@}%
	{\parindent \z@ \centering \reset@font
		\thickhrulefill\quad
		\scshape  {\chapnumfont \@chapapp{}}{\chapnumfont \thechapterwords}
		\quad \thickhrulefill
		\par\nobreak
		\vspace*{15\p@}%
		\interlinepenalty\@M
		\hrule
		\vspace*{15\p@}%
		\huge {\bfseries  #1}\par\nobreak
		\par
		\vspace*{15\p@}%
		\hrule
		\vskip 15\p@
	}}
	\def\@makeschapterhead#1{%
		\vspace*{15\p@}%
		{\parindent \z@ \centering \reset@font
			\thickhrulefill
			\par\nobreak
			\vspace*{15\p@}%
			\interlinepenalty\@M
			\hrule
			\vspace*{15\p@}%
			\Huge \bfseries #1\par\nobreak
			\par
			\vspace*{15\p@}%
			\hrule
			\vskip 30\p@
		}}
		
		\DeclareFixedFont{\chapnumfont}{T1}{phv}{b}{n}{20pt}
		\DeclareFixedFont{\chapchapfont}{T1}{phv}{b}{n}{16pt}
		\DeclareFixedFont{\chaptitfont}{T1}{phv}{b}{n}{24.88pt}
		\def\@makechapterhead#1{%
			\vspace*{15\p@}%
			{\parindent \z@ \centering \reset@font
				\thickhrulefill\quad
				\scshape {\chaptitfont\color[rgb]{0.00,0.50,1.00}\@chapapp{}}
				{\chapnumfont \thechapterwords}
				\quad \thickhrulefill
				\par\nobreak
				\vspace*{15\p@}%
				\interlinepenalty\@M
				\hrule
				\vspace*{15\p@}%
				{\Large\bfseries #1}\par\nobreak
				\par
				\vspace*{15\p@}%
				\hrule
				\vskip 30\p@
			}}%
			
\begin{document}
\title{Dynamics on Bi-Lagrangian Structures and Cherry maps}
\date{ }
\author{ 
\large{Bertuel Tangue Ndawa}\\
\normalsize{University Institute of Technology,
University of Ngaoundere},\\
\normalsize{Institut des Hautes Études Scientifiques}
\\ 
\normalsize{bertuelt@yahoo.fr
}
\vspace{0.25cm}\\
					\today}
				\maketitle
				
				\selectlanguage{english}

To the memory of the late Dr. Ben Fay\c{c}al, whose extraordinary and heroic act of selflessness at the District Hospital of Poli, in North Cameroon, led him to use nearly all of his salary arrear in order to save lives.	

\`A la m\'emoire du regrett\'e Dr. Ben Fay\c{c}al, dont l'extraordinaire et h\'ero\"{\i}que acte d'abn\'egation \`a l'H\^opital de district de Poli, dans le Nord-Cameroun, l'a conduit \`a utiliser la quasi-int\'egralit\'e de son rappel de salaire afin de sauver des vies.

\begin{center}
Part of this work was carried out during a visit to Institut des Hautes \'{E}tudes Scientifiques (IHES). We thank the scientific committee for providing this opportunity and all the IHES staff for their hospitality, in particular Fran\c{c}ois Bachelier, M\'elodie De Morais, Olivier Fernandez,
Aur\'elie Gayraud, C\'ecile Gourgues, M\'elanie Hajdinjak,
Elisabeth Jasserand, and Emmanuel Ullmo.  We also
thank  Krishnendu Gongopadhyay and his family, who were also visiting IHES, for their kindness and friendliness.
\end{center}

	\section*{Abstract}
We consider a bi-Lagrangian structure $(\omega,\mathcal{F}_1,\mathcal{F}_2)$ on a manifold $M$; that is, the quadruplet $(M,\omega,\mathcal{F}_1,\mathcal{F}_2)$ is a bi-Lagrangian manifold. We prolong bi-Lagrangian structures on $M$ and lift a given dynamics to its tangent and cotangent bundles in several ways. In some cases, we show that the lifted structures are affine. In the case of the two-dimensional torus $\mathbb{T}^2$, we prove that an extension of the same dynamics to pairs of so-called Cherry vector fields induces a conjugation action on a subset of Cherry maps, namely circle maps with a flat interval. Furthermore, we define linear connections for certain Cherry maps using bi-Lagrangian data. Finally, we show that maps generated by a bi-Lagrangian structure belong a priori to different H\"{o}lder homeomorphism classes.
\vspace{0.25cm}

				\textbf{Keywords}: Bi-Lagrangian, Cherry map, Cherry vector field, Circle map, Flat piece,
				Hess connection, H\"{o}lder,  Symplectic, Symplectomorphism.
				\vspace{0.25cm}
				
				\textbf{MSC2010}: 37C85, 37E10, 37E35, 53D05, 53D12, 93B27, 93C35, 93D40. 
				




\bigskip

\section{Introduction}
A non-degenerate and closed 2-form $\omega$ on a manifold $M$ is called a symplectic form, and the pair $(M,\omega)$ is called a symplectic manifold. A bi-Lagrangian manifold is a quadruple $(M,\omega,\mathcal{F}_1,\mathcal{F}_2)$ such that $(\mathcal{F}_1,\mathcal{F}_2)$ is a pair of transversal foliations, both Lagrangian, on the symplectic manifold $(M,\omega)$. The pair $(\mathcal{F}_1,\mathcal{F}_2)$ is then called a bi-Lagrangian structure. Note that if $(M,\omega,\mathcal{F}_1,\mathcal{F}_2)$ is a bi-Lagrangian manifold, then the triple $(\omega,\mathcal{F}_1,\mathcal{F}_2)$ is also called a bi-Lagrangian structure on $M$.

Bi-Lagrangian manifolds have been intensively studied in recent years, see \cite{GB1,1,7,TNB2}. Among the many reasons to study them are their connections to geometric quantization (see \cite{7}) and Koszul--Vinberg cohomology (see \cite{GB1}). Moreover, a bi-Lagrangian structure $(\omega,\mathcal{F}_1,\mathcal{F}_2)$ on $M$ induces, in a one-to-one manner, a para-K\"ahler structure $(G,F)$ on $M$, where $G$ is a pseudo-Riemannian metric and $F$ is the so-called almost product structure (that is, $F$ is a $(1,1)$-tensor on $M$ satisfying $F^{2}=\operatorname{Id}$ and $F\neq \operatorname{Id}$), which verifies
$G(F(\cdot),F(\cdot)) = -G(\cdot,\cdot)$.
The three tensors $\omega$, $G$, and $F$ are linked by the relation
$\omega(\cdot,\cdot) = G(F(\cdot),\cdot)$,
see \cite{1}. Therefore, a bi-Lagrangian manifold combines symplectic, semi-Riemannian, and almost product structures.

In \cite{TNB2}, the author studies a lifting of an action of the symplectic group on the set of bi-Lagrangian structures when the manifold is parallelizable. The main result is based on a lifting of bi-Lagrangian manifolds, and the action defined there is induced by the map $(\psi,X)\mapsto\psi_*X$, which is the canonical action of the diffeomorphism group on the set of vector fields. In this work, we lift the same action defined in \cite{TNB2} in a different ways, for all bi-Lagrangian manifolds.
In the case of the two-dimensional torus $\mathbb{T}^2$ (where a bi-Lagrangian structure can be defined by a pair of nowhere-vanishing vector fields that are everywhere transverse), the action $(\psi,X)\mapsto\psi_*X$ also induces a conjugation action on the set of circle maps with a flat piece (called Cherry maps), generated by certain pairs of so-called Cherry vector fields. On this basis, we define the linear connection associated with certain Cherry maps.
Note that a pair of Cherry vector fields generates several Cherry maps, and these do not belong a priori to the same H\"{o}lder class.

Before we can explain more precisely and prove our results, it is necessary to present some definitions, fix some notation   and formulate some known results we need.

\section{Notation, Basic Definitions, and Properties}\label{sub1}

We assume that all objects are smooth throughout this paper unless otherwise stated. Let $k\in\mathbb{N}^*$. Instead of $\{1,2,\dots,k\}$, we will simply write $[k]$.

Einstein summation convention: an index repeated as a subscript and superscript in a product denotes summation over the corresponding range. For example,
$$
\lambda^i\xi_i=\sum_{i=1}^n \lambda^i\xi_i,
\qquad\text{and}\qquad
X^i\partial_{y^i}=X^i\frac{\partial}{\partial y^i}=\sum_{i=1}^nX^i\frac{\partial}{\partial y^i}.
$$

\subsection{Push Forward Foliation and Conormal Distribution}
Let $M$ be an $m$-dimensional manifold. We define a $k$-dimensional foliation $\mathcal{F}$ on $M$ as a decomposition of $M$ into a union of disjoint, non-empty, connected, immersed $k$-dimensional submanifolds $\{S_x\}_{x\in M}$, called the leaves of the foliation, with the following property (the complete integrability property): every point $y\in M$ has a coordinate chart $(U, y^1,\dots,y^m)$ such that, for each leaf $S_x$, the components of $U\cap S_x$ are described by the equations
$y^{k+1}=\text{constant},\quad \dots,\quad y^m=\text{constant}$.

Let $TM$ be the tangent bundle of $M$, and let $\mathfrak{X}(M)$ denote the set of vector fields on $M$, i.e., $\mathfrak{X}(M)=\Gamma(TM)$. The expressions $T\mathcal{F}\subset TM$ and $\Gamma(\mathcal{F})=\Gamma(T\mathcal{F})\subset \mathfrak{X}(M)$ denote the tangent bundle of $\mathcal{F}$ and the space of sections of $T\mathcal{F}$, respectively. For each point $y\in M$, the vector subspace $T_yS_y\subset T_yM$ is called the tangent bundle of $\mathcal{F}$ at $y$, and is denoted by $\mathcal{F}_y$ (or $T_y\mathcal{F}$).

The Lie bracket of two vector fields $X,Y$ on $M$ is given by
$[X,Y] := X\circ Y - Y\circ X$.
The complete integrability property of a foliation $\mathcal{F}$ means that $\Gamma(\mathcal{F})$ is stable under the Lie bracket, i.e., if $X,Y\in\Gamma(\mathcal{F})$, then $[X,Y]\in\Gamma(\mathcal{F})$. This is the Frobenius integrability theorem.

Let $\psi:M\to N$ be a diffeomorphism. The pushforward foliation
$\psi_*\mathcal{F}=\{\psi(S_x)\}_{x\in M}$ is a foliation on $N$, and
\begin{equation}
	\Gamma(\psi_*\mathcal{F}) := \{\psi_*X,\; X\in \Gamma(\mathcal{F})\} = \psi_*\Gamma(\mathcal{F}).\label{Bieq2}
\end{equation}
Let $S\subset M$ be a submanifold. The conormal space at a point $x\in S$ is defined by
$N_x^*S=\{\xi_x\in T_xM:\ \xi_x|_{T_xS}=0\}$,
and the conormal bundle of $S$ is
$N^*S=\{(\xi_x,x)\in T^*M:\ \xi_x\in N_x^*S\}$.
The conormal bundle of a foliation $\mathcal{F}=\{S_x\}_{x\in M}$ is
$N^*\mathcal{F}=\{N^*S_x\}_{x\in M}$.
For more details, see \cite[p.\,17--18]{dasilva}.

\subsection{Bi-Lagrangian Manifold}
Let $(M,\omega)$ be a $2n$-dimensional symplectic manifold. Since $\omega$ is non-degenerate, it defines a vector bundle isomorphism
$\omega^\flat:TM\longrightarrow T^*M$, $(x,X_x)\longmapsto (x,\omega_x(X_x,\cdot))$,
and we denote its inverse by $\omega^\sharp:T^*M\longrightarrow TM$.

A foliation $\mathcal{F}$ is Lagrangian if its orthogonal section
$$\bigl(\Gamma(T\mathcal{F})\bigr)^\perp:=\left\{X\in\mathfrak{X}(M):\ \omega(X,Y)=0,\ \forall Y\in\Gamma(T\mathcal{F})\right\}$$
is equal to $\Gamma(T\mathcal{F})$.

A bi-Lagrangian structure on $M$ consists of a pair $(\mathcal{F}_1,\mathcal{F}_2)$ of transversal Lagrangian foliations, together with a symplectic form $\omega$. As a consequence, the tangent bundle decomposes as
$TM=T\mathcal{F}_1\oplus T\mathcal{F}_2$.

Let $(\mathcal{F}_1,\mathcal{F}_2)$ be a bi-Lagrangian structure on a $2n$-dimensional symplectic manifold $(M,\omega)$. Every point in $M$ has an open neighborhood $U$ that is the domain of a chart with local coordinates $(p^1,\dots,p^n,q^1,\dots,q^n)$ such that
$$\Gamma(\mathcal{F}_1)_{\mid U}=\left<\frac{\partial}{\partial p^1},\dots,\frac{\partial}{\partial p^n}\right>
\quad \text{and}\quad
\Gamma(\mathcal{F}_2)_{\mid U}=\left<\frac{\partial}{\partial q^1},\dots,\frac{\partial}{\partial q^n}\right>.$$
Such a chart is said to be adapted to the bi-Lagrangian structure $(\mathcal{F}_1,\mathcal{F}_2)$. Moreover, if
$\omega=\sum_{i=1}^{n}dq^i\wedge dp^i$,
then such a chart is said to be adapted to the bi-Lagrangian structure $(\omega,\mathcal{F}_1,\mathcal{F}_2)$.

\subsection{Linear Connection}
Let $\nabla$ be a linear connection. The torsion tensor $T_{\nabla}$ (or simply $T$ if there is no ambiguity) and the curvature tensor $R_{\nabla}$ (or simply $R$) are defined, respectively, by
\begin{equation*}
	T_{\nabla}(X,Y)=\nabla_XY-\nabla_YX-[X,Y], \quad X,Y\in\mathfrak{X}(M),
\end{equation*}
and
\begin{equation*}
	R_{\nabla}(X,Y)Z=\nabla_X\nabla_YZ-\nabla_Y\nabla_XZ-\nabla_{[X,Y]}Z,\quad 
	X,Y,Z\in\mathfrak{X}(M).\label{courbure}
\end{equation*}

We say that a linear connection $\nabla$
\begin{enumerate}
	\item[-] parallelizes $\omega$ if $\nabla\omega=0$, that is,
	$\omega(\nabla_XY,Z)+\omega(Y,\nabla_XZ)=X\bigl(\omega(Y,Z)\bigr)$
	for every $X,Y,Z\in\mathfrak{X}(M)$;
	\item[-] preserves $\mathcal{F}$ if $\nabla\Gamma(\mathcal{F})\subseteq \Gamma(\mathcal{F})$, meaning that
	$\nabla_XY\in\Gamma(\mathcal{F})$ for every $X\in\mathfrak{X}(M)$ and $Y\in\Gamma(\mathcal{F})$.
\end{enumerate}


\subsection{Prolongations of Geometric Objects to the Tangent Bundle}
We recall here the essential results on complete and vertical lifts. For more
details, the reader is referred to \cite{YKI,YKII}. If $\zeta$ is an object on a manifold $M$, then $\zeta^v$ and $\zeta^c$ denote, respectively, the vertical and complete lifts of $\zeta$ to $TM$.

Let $\pi:TM\to M$ be the natural projection. For $f\in C^\infty(M)$, $X\in\mathfrak{X}(M)$, and $\alpha\in\Omega^1(M)$, their vertical and complete lifts to $TM$ are defined as follows:
\begin{enumerate}
	\item[-] $f^v=f\circ\pi$, and $f^c=df$, where $df$ is regarded as a fiberwise linear function on $TM$;
	\item[-] $X^v(f^c)=(Xf)^v$, and $X^c(f^c)=(Xf)^c$;
	\item[-] $\alpha^v(Y^c)=(\alpha(Y))^v$, and $\alpha^c(Y^c)=(\alpha(Y))^c$ for every $Y\in\mathfrak{X}(M)$.
\end{enumerate}

From these definitions, it follows that, for every $f,g\in C^\infty(M)$, $X,Y\in\mathfrak{X}(M)$, and $\alpha\in\Omega^1(M)$,
\begin{equation*}\label{eqlift1}
	\begin{array}{lll}
		(fg)^v=f^v g^v, &(f\alpha)^v=f^v\alpha^v, &
		(fX)^v=f^v X^v, \\[2mm] (fX)^c=f^c X^v+f^v X^c,&
		(fg)^c=f^c g^v+f^v g^c,&
		(f\alpha)^c=f^c \alpha^v+f^v \alpha^c.
	\end{array}
\end{equation*}
Moreover,
\begin{equation}\label{eqlift4}
	[X^c,Y^c]=[X,Y]^c,\qquad
	[X^v,Y^v]=0,\qquad
	[X^c,Y^v]=[X,Y]^v.
\end{equation}
In local coordinates $(x^1,\ldots,x^n,y^1,\ldots,y^n)$ on $TM$, if $X=X^i\frac{\partial}{\partial x^i}$, then
\begin{equation}\label{eqlift5}
	X^c=
	X^i\frac{\partial}{\partial x^i}
	+\frac{\partial X^i}{\partial x^j}\,y^j\frac{\partial}{\partial y^i},
	\qquad
	X^v=
	X^i\frac{\partial}{\partial y^i}.
\end{equation}

Let $S$ and $T$ be two tensor fields on $M$. The vertical and complete lifts of tensor fields are defined inductively by the following formulas:
\begin{equation*}\label{eqlift2}
	\begin{array}{ll}
		(S+T)^v=S^v+T^v, &
		(S+T)^c=S^c+T^c,\\[2mm]
		(S\otimes T)^v=S^v\otimes T^v, &
		(S\otimes T)^c=S^c\otimes T^v+S^v\otimes T^c.
	\end{array}
\end{equation*}

\begin{remark}\label{DPDrem1}
	Let $\mathcal{F}$ be a foliation on $M$, and denote by the same symbol
	$\mathcal{F}$ its tangent distribution. Since $\mathcal{F}$ is involutive,
	it follows from the bracket relations \eqref{eqlift4} that the distribution
	on $TM$ defined by
	\begin{equation}\label{lift of foliation}
		\Gamma(\mathcal{F}^{cv})
		=
		\left\langle X^c,\;X^v;\;X\in\Gamma(\mathcal{F})\right\rangle_{C^\infty(TM)}
	\end{equation}
	is involutive. Hence, by Frobenius' theorem, it induces a foliation
	$\mathcal{F}^{cv}$ on $TM$, called the total lift of $\mathcal{F}$.
	
	Moreover, using the local expressions in \eqref{eqlift5}, we obtain
	\begin{align}
		\Gamma(\mathcal{F}^{cv})
		&=
		\left\langle X^c;\;X\in\Gamma(\mathcal{F})\right\rangle_{C^\infty(TM)}
		\oplus
		\left\langle X^v;\;X\in\Gamma(\mathcal{F})\right\rangle_{C^\infty(TM)}
		\nonumber\\
		&=
		\Gamma(\mathcal{F}^c)\oplus\Gamma(\mathcal{F}^v).
		\label{eqlift6}
	\end{align}
	In particular,
	$$
	\operatorname{rank}(\mathcal{F}^c)
	=
	\operatorname{rank}(\mathcal{F}^v)
	=
	\operatorname{rank}(\mathcal{F}).
	$$
	Combining this with \eqref{eqlift6}, we obtain
	\begin{equation}\label{eqlift7}
		\operatorname{rank}(\mathcal{F}^{cv})
		=
		2\operatorname{rank}(\mathcal{F}^c)
		=
		2\operatorname{rank}(\mathcal{F}^v).
	\end{equation}
\end{remark}






\section{Technical Tools}
\subsection{Symplectic Manifolds}
The cotangent bundle $T^*M$ carries a canonical symplectic structure given by the tautological 2-form. This section states the result precisely, see \cite{dasilva} for more details.

Let $M$ be an $m$-dimensional manifold, and let ${}^*\hspace{-0.1cm}\pi:T^*M\longrightarrow M$ be the natural projection.
The tautological 1-form, also known as the Liouville 1-form, $\theta$, is defined by
$$
\theta_{(x,\alpha_x)}(v)=\alpha_x\!\left((T_{(x,\alpha_x)}{}^*\hspace{-0.1cm}\pi)(v)\right),
\qquad (x,\alpha_x)\in T^*M,\ v\in T_{(x,\alpha_x)}T^*M.
$$
The exterior differential $d\theta$ of $\theta$ is called the canonical symplectic form, also known as the Liouville 2-form, on the cotangent bundle $T^*M$.

Note that for any coordinate chart $(U,x^1,\dots,x^m)$ on $M$, with the associated cotangent coordinate chart
$(T^*U,x^1,\dots,x^m,\xi_1,\dots,\xi_m)$, we have
$$
\theta=\sum_{i=1}^{m}\xi_i\,dx^i
\qquad\text{and}\qquad
d\theta=\sum_{i=1}^{m}d\xi_i\wedge dx^i.
$$

\begin{proposition}\label{lifting of symplecto}
	Let $\psi:M_1\longrightarrow M_2$ be a diffeomorphism. The lift
	$\hat{\psi}:T^*M_1\to T^*M_2$, defined by 	$
	(x,\alpha_x)\longmapsto\bigl(\psi(x),(\psi^{-1})^*\alpha_x\bigr),
	$
	is a symplectomorphism from $(T^*M_1,d\theta_1)$ to $(T^*M_2,d\theta_2)$, where $d\theta_1$ and $d\theta_2$ are the canonical symplectic forms on $T^*M_1$ and $T^*M_2$, respectively.
\end{proposition}

\subsection{Affine and Push Forward Bi-Lagrangian Structures}
In \cite{7}, the author shows that every bi-Lagrangian manifold admits a unique torsion-free connection, called the Hess connection or bi-Lagrangian connection, which parallelizes $\omega$ and preserves both foliations. Moreover, the author provides an explicit formula for this connection.

A bi-Lagrangian manifold is said to be affine if the curvature tensor of its Hess connection vanishes identically. This condition is necessary and sufficient for defining the Koszul--Vinberg cohomology, see \cite{GB1}.
Affine bi-Lagrangian manifolds are characterized as follows:

\begin{theorem}\cite[Theor.~2, p.~159]{7}\label{c10}
	Let $(\omega,\mathcal{F}_1,\mathcal{F}_2)$ be a bi-Lagrangian structure on a $2n$-dimensional manifold $M$, and let $\nabla$ denote its Hess connection. The following assertions are equivalent:
	\begin{enumerate}
		\item[(a)] The connection $\nabla$ is flat;
		\item[(b)] Each point of $M$ has a coordinate chart adapted to $(\omega,\mathcal{F}_1,\mathcal{F}_2)$.
	\end{enumerate}
\end{theorem}

The following result describes a method for pushing forward a bi-Lagrangian structure, illustrating how these structures transform under a symplectomorphism.

\begin{lemma}\cite[Lem.~2.2, p.~6]{TNB2}\label{Bilem1}
	Let $(M,\omega,\mathcal{F}_1,\mathcal{F}_2)$ be a bi-Lagrangian manifold with $\nabla$ as its Hess connection, and let $N$ be a manifold diffeomorphic to $M$. Then for any diffeomorphism $\psi:M\longrightarrow N$, the structure 	$
	\bigl((\psi^{-1})^*\omega,\ \psi_*\mathcal{F}_1,\ \psi_*\mathcal{F}_2\bigr)
	$
	is bi-Lagrangian on $N$, with Hess connection
	$
	\nabla^{\psi}:(X,Y)\longmapsto \psi_*\bigl(\nabla_{\psi_*^{-1}X}(\psi_*^{-1}Y)\bigr).
	$
	Moreover, if $(\omega,\mathcal{F}_1,\mathcal{F}_2)$ is affine, then so is
	$
	\bigl((\psi^{-1})^*\omega,\ \psi_*\mathcal{F}_1,\ \psi_*\mathcal{F}_2\bigr).
	$
\end{lemma}




\section{Statements and Proofs of Results}

\subsection{Lift of a Dynamic Defined on the Set Bi-Lagrangian Structures}
Our first result provides a construction of lifted bi-Lagrangian structures on the tangent bundle $TM$ and the cotangent bundle $T^*M$.

\begin{theorem}\label{Bitheo1}
	Let $(M,\omega,\mathcal{F}_1,\mathcal{F}_2)$ be a bi-Lagrangian manifold. We have:
	\begin{enumerate}
		\item the quadruplet $(T^*M,d\theta,N^*\mathcal{F}_1,N^*\mathcal{F}_2)$ is an affine bi-Lagrangian manifold;
		\item the expression $(T^*M,\tilde{\omega}:={}^{*}\hspace{-0.1cm}\pi^*\omega+d\theta,N^*\mathcal{F}_1,N^*\mathcal{F}_2)$ is a bi-Lagrangian manifold;
		\item if $(M,\omega,\mathcal{F}_1,\mathcal{F}_2)$ is affine, then $(TM,(\omega^\flat)^*d\theta,\mathcal{F}_1^{cv},\mathcal{F}_2^{cv})$ is an affine bi-Lagrangian manifold.
	\end{enumerate}
\end{theorem}

\begin{remark}\label{Birem4}
	Let $(M,\omega,\mathcal{F}_1,\mathcal{F}_2)$ be a bi-Lagrangian manifold.
	\begin{enumerate}
		\item Since $TM$ is diffeomorphic to $T^*M$ via $\omega^\flat$, then by Lemma~\ref{Bilem1} a bi-Lagrangian structure on $TM$ induces a bi-Lagrangian structure on $T^*M$, and vice versa.
		\item Let $\tilde{L}_i(\psi)=\hat{\psi}$ for $i=1,2$, and let $\tilde{L}_3(\psi)=\omega^{\#}\circ \hat{\psi}\circ \omega^{\flat}$. Moreover,
		\[
		\tilde{L}_i(M,\omega,\mathcal{F}_1,\mathcal{F}_2)=
		\begin{cases}
			(T^*M,d\theta,N^*\mathcal{F}_1,N^*\mathcal{F}_2) & \text{if }\quad i=1,\\
			(T^*M,\tilde{\omega},N^*\mathcal{F}_1,N^*\mathcal{F}_2) & \text{if }\quad  i=2,\\
			(TM,(\omega^\flat)^*d\theta,\mathcal{F}_1^{cv},\mathcal{F}_2^{cv}) & \text{if }\quad i=3.
		\end{cases}
		\]
		Let $(\psi,(\mathcal{F}_1,\mathcal{F}_2))\mapsto(\psi_*\mathcal{F}_1,\psi_*\mathcal{F}_2)$ be the action of the symplectic group on the set of bi-Lagrangian structures on $M$, see \cite[Theo.~2.3, p.~414]{TNB2}. The following holds:
		for every $i\in\{1,2,3\}$, if $(\tilde{L}_i(\psi))_*\mathcal{F}_j\subset \tilde{L}_i(\psi_*\mathcal{F}_j)$ for some $j\in\{1,2\}$, the following diagram is commutative.
		
		\vspace{0.75cm}
		\begin{center}
			\setlength{\unitlength}{1mm}
			\thicklines
			\begin{picture}(40,20)
				\put(-20,0){$(M,\omega,\mathcal{F}_1,\mathcal{F}_2)$}
				\put(36.5,0){$(M,\omega,\psi_*\mathcal{F}_1,\psi_*\mathcal{F}_2)$}
				\put(-23,20){$\tilde{L}_i(M,\omega,\mathcal{F}_1,\mathcal{F}_2)$}
				\put(31,20){$\tilde{L}_i(M,\omega,\psi_*\mathcal{F}_1,\psi_*\mathcal{F}_2)$}
				\put(18,-2){$\psi_*$}
				\put(11,23){$(\tilde{L}_i(\psi))_*$}
				\put(-13.5,7.5){\rotatebox{90}{Lift}}
				\put(7.5,3){Push forward}
				\put(-10,3.5){\vector(0,1){15.5}}
				\put(48,3.5){\vector(0,1){15.5}}
				\put(4.5,21.55){\vector(1,0){26.5}}
				\put(3.5,1.55){\vector(1,0){32.5}}
			\end{picture}
		\end{center}
		Moreover, by Lemma~\ref{Bilem1}, Theorem~\ref{Bitheo1}, and
		Remark~\ref{DPDrem1}, for $i=1,3$, if the bi-Lagrangian structure
		$(\omega,\mathcal{F}_1,\mathcal{F}_2)$ is affine, then the three
		bi-Lagrangian structures induced in the diagram are also affine.
	\end{enumerate}
\end{remark}





\begin{proof}[Proof of Theorem~\ref{Bitheo1}]
	The following result will be useful in the proof.
	\begin{lemma}\label{Bilem2}
		Let $\mathcal{F}=\{S_x\}_{x\in M}$ be a $k$-foliation on an $m$-dimensional manifold $M$. Then $N^*\mathcal{F}$ is a Lagrangian foliation on $(T^*M,d\theta)$. Moreover, if $M$ is endowed with a symplectic form $\omega$, and $\mathcal{F}$ is Lagrangian on $(M,\omega)$, then $N^*\mathcal{F}$ is Lagrangian on $(T^*M,\tilde{\omega})$.
	\end{lemma}
	
	\begin{proof}
		First, note that if $S\subset M$ is a submanifold, then $N^*S$ is a Lagrangian submanifold of $(T^*M,d\theta)$ (this is a general result, see \cite[Cor.~3.7, p.~18]{dasilva}). As a consequence, $N^*\mathcal{F}=\{N^*S_x\}_{x\in M}$ is Lagrangian on $(T^*M,d\theta)$.
		Hence it remains to show that $N^*\mathcal{F}$ is completely integrable.
		
		Since $N^*\mathcal{F}$ is Lagrangian, complete integrability is equivalent to
		\begin{equation*}
			d\theta([X,Y],Z)=0
			\qquad\text{for all }X,Y,Z\in\Gamma(N^*\mathcal{F}).
		\end{equation*}
		Let $(U,p^1,\dots,p^m)$ be a coordinate chart adapted to the foliation $\mathcal{F}$, and let $$(T^*U,p^1,\dots,p^m,\xi_1,\dots,\xi_m)$$ be the associated coordinate chart on $T^*M$. Observe that
		\begin{equation}\label{Bieq8}
			\Gamma(N^*\mathcal{F})\big|_{T^*U}
			=\left\langle
			\frac{\partial}{\partial p^1},\dots,\frac{\partial}{\partial p^k},
			\frac{\partial}{\partial \xi_{k+1}},\dots,\frac{\partial}{\partial \xi_m}
			\right\rangle.
		\end{equation}
		Let us write
		\begin{equation*}
			\begin{cases}
				(y^i)_{i=1,\dots,m}=\bigl((p^i)_{i=1,\dots,m},(\xi_i)_{i=k+1,\dots,m}\bigr),\\
				X=X^i\frac{\partial}{\partial y^i},\qquad
				Y=Y^j\frac{\partial}{\partial y^j},\qquad
				Z=Z^k\frac{\partial}{\partial y^k}.
			\end{cases}
		\end{equation*}
		Then
		\begin{equation*}
			\begin{cases}
				[X,Y]=\mu^j\frac{\partial}{\partial y^j},\\
				[[X,Y],Z]=\lambda^j\frac{\partial}{\partial y^j},
			\end{cases}
			\quad \text{ where }\quad
			\begin{cases}
				\mu^j=X^i\frac{\partial Y^j}{\partial y^i}-Y^i\frac{\partial X^j}{\partial y^i},\\
				\lambda^j=\mu^i\frac{\partial Z^j}{\partial y^i}-Z^i\frac{\partial \mu^j}{\partial y^i}.
			\end{cases}
		\end{equation*}
		Therefore,
		\begin{align*}
			[X,Y]\theta(Z)&=\mu^i\frac{\partial}{\partial
				y^i}(Z^k\xi_k),\\
			\theta([[X,Y],Z])&=\lambda^j\xi_j, 
			\\Z\theta([X,Y])&=Z^k\frac{\partial}{\partial y^k}(\mu^i\xi_i). 
		\end{align*}
		Hence,
		\begin{equation*}
			d\theta([X,Y],Z)
			= [X,Y]\theta(Z)-\theta([[X,Y],Z])-Z\theta([X,Y])
			=0.
		\end{equation*}
		So $N^*\mathcal{F}$ is a Lagrangian foliation on $(T^*M,d\theta)$.
		
		Now, suppose that $\mathcal{F}$ is a Lagrangian foliation on $(M,\omega)$.
		Observe that, if $X\in\Gamma(N^*\mathcal{F})$, then
		${}^*\pi_*X\in\Gamma(\mathcal{F})$. Hence, for every
		$X,Y\in\Gamma(N^*\mathcal{F})$, we have 
		\begin{equation*}
			\tilde{\omega}(X,Y)=\omega(^*\hspace{-0.1cm}\pi_*X,^*\hspace{-0.1cm}\pi_*Y )+d\theta(X,Y)=0
		\end{equation*}
		where we have used the fact that $\mathcal{F}$ is Lagrangian, together with the
		definition of $d\theta$. Therefore, $N^*\mathcal{F}$ is a Lagrangian foliation
		on $(T^*M,\widetilde{\omega})$. This completes the proof of
		Lemma~\ref{Bilem2}.
	\end{proof}

	We now prove Theorem~\ref{Bitheo1}. Let $(M,\omega,\mathcal{F}_1,\mathcal{F}_2)$ be a bi-Lagrangian manifold.
	\begin{enumerate}
		\item By Lemma~\ref{Bilem2}, $(N^*\mathcal{F}_1,N^*\mathcal{F}_2)$ is a pair of Lagrangian foliations on $(T^*M,d\theta)$.
		From \eqref{Bieq8}, it follows that if $(U,p^1,\dots,p^{n},q^1,\dots,q^{n})$ is a coordinate chart adapted to the bi-Lagrangian structure $(\mathcal{F}_1,\mathcal{F}_2)$, with $(T^*U,p^1,\dots,p^{n},q^1,\dots,q^{n},\xi_1,\dots,\xi_{2n})$ as the corresponding bundle coordinate chart on $T^*M$, then
		\begin{equation}\label{Bieq9}
			\begin{cases}
				\Gamma(N^*\mathcal{F}_1)\big|_{T^*U}
				=\left<\frac{\partial}{\partial p^1},\dots,\frac{\partial}{\partial p^n},
				\frac{\partial}{\partial \xi_{n+1}},\dots,\frac{\partial}{\partial \xi_{2n}}\right>,\vspace{0.25cm}\\
				\Gamma(N^*\mathcal{F}_2)\big|_{T^*U}
				=\left<\frac{\partial}{\partial q^{1}},\dots,\frac{\partial}{\partial q^{n}},
				\frac{\partial}{\partial \xi_1},\dots,\frac{\partial}{\partial \xi_n}\right>,
			\end{cases}
		\end{equation}
		and
		\begin{equation}\label{Bieq10}
			d\theta\big|_{T^*U}
			=\sum_{i=1}^{n}\left(d\xi_i\wedge dp^i+d\xi_{n+i}\wedge dq^i\right).
		\end{equation}
		
		Therefore, by \eqref{Bieq9}, the foliations $N^*\mathcal{F}_1$ and $N^*\mathcal{F}_2$ are transverse.
		Since they are Lagrangian with respect to $d\theta$ by Lemma~\ref{Bilem2}, the triple
		$(d\theta,N^*\mathcal{F}_1,N^*\mathcal{F}_2)$
		is a bi-Lagrangian structure on $T^*M$.
		Finally, combining \eqref{Bieq9} and \eqref{Bieq10}, Theorem~\ref{c10} implies that
		$(d\theta,N^*\mathcal{F}_1,N^*\mathcal{F}_2)$ is affine.
		
		\item The second assertion of Theorem~\ref{Bitheo1} follows directly by combining Lemma~\ref{Bilem2} with the equalities in \eqref{Bieq9}.
		
		\item Now suppose that $(M,\omega,\mathcal{F}_1,\mathcal{F}_2)$ is affine.
		Let $(p^1,\dots,p^{n},q^1,\dots,q^{n})$ be a coordinate chart adapted to $(\omega,\mathcal{F}_1,\mathcal{F}_2)$, and let $$(p^1,\dots,p^{n},q^1,\dots,q^{n},x^1,\dots,x^{2n})$$ be its associated tangent bundle coordinates.
		We have
		\begin{equation}\label{eqlift8}
			\Gamma(\mathcal{F}_1)\big|_{U}
			=\left<\frac{\partial}{\partial p^1},\dots,\frac{\partial}{\partial p^n}\right>,
			\qquad
			\Gamma(\mathcal{F}_2)\big|_{U}
			=\left<\frac{\partial}{\partial q^{1}},\dots,\frac{\partial}{\partial q^{n}}\right>
		\end{equation}
		and
		\begin{equation}\label{eqlift9}
			\left(\omega^\flat\right)^*(d\theta)
			=\sum_{i=1}^n dx^{n+i}\wedge dp^i-\sum_{i=1}^n dx^{i}\wedge dq^i.
		\end{equation}
		
		Combining \eqref{eqlift5}, \eqref{lift of foliation} (or \eqref{eqlift6}), \eqref{eqlift7}, \eqref{eqlift8}, and \eqref{eqlift9}, we conclude that $\mathcal{F}_i^{cv}$, $i\in[2]$, are Lagrangian distributions.
		By \eqref{eqlift4} and the Frobenius theorem, they are foliations.
		
		Moreover, since the coefficients of $\left(\omega^\flat\right)^*(d\theta)$ are constant, see \eqref{eqlift9}, the Christoffel symbols of the Hess connection associated with $$(TM,\left(\omega^\flat\right)^*d\theta,\mathcal{F}_1^{cv},\mathcal{F}_2^{cv})$$ vanish.
		Indeed, this connection is torsion-free, parallelizes $\left(\omega^\flat\right)^*d\theta$, preserves both foliations, and, by uniqueness, it is the Hess connection of the structure.
		This completes the proof of the theorem.
	\end{enumerate}
\end{proof}

\begin{example}
	Let $M=\mathbb{R}^2$ and $\omega=dy\wedge dx$. Let
	$$
	\Gamma(\mathcal{F}_1)=\left\langle \frac{\partial}{\partial x}\right\rangle,
	\qquad
	\Gamma(\mathcal{F}_2)=\left\langle \frac{\partial}{\partial y}\right\rangle.
	$$
	Let $(x,y)$, $(x,y,u,w)$, and $(x,y,z,t)$ be the coordinates on $M$, $T^*M$, and $T^*M$, respectively. We have
	$$
	d\theta = dz \wedge dx + dt \wedge dy,\qquad
	^*\hspace{-0.1cm}\pi^*\omega = dy \wedge dx,\qquad
	\left(\omega^\flat\right)^*(d\theta) = dw \wedge dx - du \wedge dy.
	$$
	Moreover,
	$$
	\Gamma(N^*\mathcal{F}_1)=\left\langle \frac{\partial}{\partial x}, \frac{\partial}{\partial t} \right\rangle,
	\qquad
	\Gamma(N^*\mathcal{F}_2)=\left\langle \frac{\partial}{\partial y}, \frac{\partial}{\partial z} \right\rangle,
	$$
	and
	$$
	\Gamma(\mathcal{F}_1^{cv}) = \left\langle \frac{\partial}{\partial x}, \frac{\partial}{\partial u} \right\rangle,
	\qquad
	\Gamma(\mathcal{F}_2^{cv}) = \left\langle \frac{\partial}{\partial y}, \frac{\partial}{\partial w} \right\rangle.
	$$
\end{example}

\subsection{2-Dimensional Torus: Linear Connection of a Cherry  Map}
Since $\mathbb{T}^2$ is a two-dimensional orientable manifold, it is endowed with a symplectic form $\omega$.
A bi-Lagrangian foliation on $(\mathbb{T}^2,\omega)$ can be represented by two transverse nonvanishing vector fields $(X,Y)$.
The push-forward action of diffeomorphisms,
$(\psi,(X,Y))\mapsto(\psi_*X,\psi_*Y)$,
therefore acts on bi-Lagrangian structures.
Restricting this action to Cherry vector fields yields a conjugation action on certain circle maps with a flat interval, possibly having different critical exponents.

Let $X$ be a vector field on $\mathbb{T}^2$, and let $x\in \mathbb{T}^2$.
We use the following notation:
\begin{enumerate}
	\item[-] $\operatorname{Sing}(X)$ denotes the set of singularities of $X$.
	\item[-] $t\mapsto \Phi_X^t$ denotes the flow of $X$.
	Thus, for every $y$ in the domain of the flow,
	$$
	\frac{d}{dt}\Phi_X^t(y)=X_{\Phi_X^t(y)}.
	$$
	\item[-] 
	$\gamma^-(x)=\{\Phi_X^t(x): t\leq 0\}$ and
	$\gamma^+(x)=\{\Phi_X^t(x): t\geq 0\}$
	denote, respectively, the negative and positive semi-trajectories of $x$.
	\item[-] $\gamma(x)=\gamma^-(x)\cup \gamma^+(x)$
	denotes the trajectory of $x$.
\end{enumerate}

A trajectory $\gamma(x)$ is said to be closed, or periodic, if it is
homeomorphic to the circle $\mathbb{S}^1$. Equivalently, $\gamma(x)$ is
periodic if there exists $T>0$ such that $\Phi_X^T(x)=x$ and $x$ is not a
singular point of $X$. A trajectory $\gamma(x)$ is said to be
non-closed if it is neither a fixed point nor a periodic trajectory.

A vector field on $\mathbb{T}^2$ without periodic trajectories is called a
Cherry vector field if it has exactly two singularities: one sink and one
saddle, both of which are hyperbolic. We denote the set of such vector fields
by $\mathfrak{X}_c(\mathbb{T}^2)$. The first example of a Cherry vector field
was constructed by Cherry, see \cite{TC}.

In this work, we relate certain pairs of Cherry vector fields to a subclass
of $\mathscr{L}$ consisting of circle maps with a flat piece. Identifying
$\mathbb{S}^1$ with $[0,1]/\{0\sim1\}$, let $\ell_1,\ell_2\geq0$ and
$U=(a,b)\subset[0,1]$. An orientation-preserving map $f\in\mathscr{L}$
has critical exponents $(\ell_1,\ell_2)$ and flat piece $U$ if the following
conditions are satisfied:
\begin{enumerate}
	\item The image of $U$ is one point.
	\item The restriction of $f$ to $[0,1]\setminus\overline{U}$ is a
	diffeomorphism onto its image.
	\item On a left-sided neighborhood of $a$, $f$ equals
	$h_l((x-a)^{\ell_1})$, where $h_l$ is a diffeomorphism on a two-sided
	neighborhood of $a$.
	\item Analogously, on some right-sided neighborhood of $b$, $f$ can be
	represented as $h_r((x-b)^{\ell_2})$.
\end{enumerate}

A map in $\mathscr{L}$ is called a Cherry map. We say that a map $f\in \mathscr{L}$
with critical exponents $(\ell_1,\ell_2)$ is symmetric if $\ell_1=\ell_2$.

Let $X$ be a Cherry vector field on
$\mathbb{T}^2\simeq \mathbb{S}^1\times[0,1]/\sim$, where
$\mathbb{S}^1\times\{0\}$ is identified with $\mathbb{S}^1\times\{1\}$.
Let $\Delta\subset \mathbb{S}^1\times\{0\}$ be the set of points whose positive
orbit reaches $\mathbb{S}^1\times\{1\}$, and let $t(x)>0$ be the first hitting
time. Then there exists $c\in\mathbb{S}^1$ such that the map
\begin{equation*}
	f:\mathbb{S}^1 \longrightarrow \mathbb{S}^1,\quad
	x\longmapsto
	\begin{cases}
		\Phi_X^{t(x)}(x) & \text{if } x\in \Delta,\\
		c & \text{if } x\notin \Delta,
	\end{cases}
\end{equation*}
is a symmetric map in $\mathscr{L}$, see \cite[\S 6.2, 6.3, p. 148-149]{LP1}.

\begin{theorem}~\label{Bitheo2}
	\begin{enumerate}
		\item Every pair $(X_1,X_2)$ of Cherry vector fields generates symmetric and non-symmetric maps $f_{X_1,X_2}$ in $\mathscr{L}$.
		The maps generated by $(X_1,X_2)$ do not necessarily belong, a priori, to the same H\"{o}lder homeomorphism class.
		
		\item Let $\mathfrak{X}_c^2(\mathbb{T}^2)=\mathfrak{X}_c(\mathbb{T}^2)\times\mathfrak{X}_c(\mathbb{T}^2)$.
		The push-forward of vector fields defines a left action:
		\begin{equation}\label{Biact1}
			*: \operatorname{Diff}(\mathbb{T}^2)\times\mathfrak{X}_c^2(\mathbb{T}^2)\longrightarrow \mathfrak{X}_c^2(\mathbb{T}^2),\quad
			(\psi, (X_1,X_2))\mapsto (\psi_*X_1,\psi_*X_2).
		\end{equation}
		This action induces a conjugation action on the corresponding class of Cherry maps:
		\begin{equation}\label{Biact2}
			\circ:\operatorname{Diff}(\mathbb{S}^1)\times\mathscr{L}_c\longrightarrow\mathscr{L}_c,\quad
			(\varphi, f)\mapsto \varphi\circ f\circ\varphi^{-1},
		\end{equation}
		where $\mathscr{L}_c$ denotes the subset of $\mathscr{L}$ consisting of maps generated by pairs of Cherry vector fields.
	\end{enumerate}
\end{theorem}

\begin{proof}
	\begin{enumerate}
		\item Let $X_1$ and $X_2$ be two Cherry vector fields. For each $i=1,2$, the
		vector field $X_i$ generates a map $f_i\in\mathscr{L}$ with flat piece
		$U_i=(a_i,b_i)$ and critical exponents $(\ell_i,\ell_i)$. The maps
		$\left(f_i\circ f_j\right)^{\ell_l,\ell_d}$, for $i,j\in[2]$, are circle maps with a flat
		piece and critical exponents $(\ell_l,\ell_d)$ of the form
		$(\ell_{i_1},\ell_{i_2})$, $(\ell_{i_1},\ell_{1}\ell_{2})$, or
		$(\ell_{1}\ell_{2},\ell_{i_2})$, where $\ell_{i_1},\ell_{i_2}\in\{\ell_1,\ell_2\}$,
		with all possible choices.
		
		We may assume, via a translation and the conjugation action on the maps,
		that $f_1(U_1)=f_2(U_2)$ and that $U:=U_1\cup U_2$ is an interval. Thus,
		for $i,j\in[2]$, the map
		\begin{equation*}
			f_{ij}(x)=
			\begin{cases}
				f_i(x) & \text{if } x \text{ lies in the left component of } [0,1]\setminus U,\\
				f_1(U_1) & \text{if } x\in U,\\
				f_j(x) & \text{if } x \text{ lies in the right component of } [0,1]\setminus U
			\end{cases}
		\end{equation*}
		belongs to $\mathscr{L}$, with flat piece $U$ and critical exponents $(\ell_i,\ell_j)$.
		
		Assume now that
		\begin{equation*}
			\left(f_i\circ f_j\right)^{\ell_{1e},\ell_{1d}}
			\quad\text{and}\quad
			\left(f_i\circ f_j\right)^{\ell_{2e},\ell_{2d}}
		\end{equation*}
		have the same rotation number, for instance, the golden mean rotation number;
		see Example~\ref{Biex1}. Define
		\begin{equation*}
			\lambda(\ell_1,\ell_2)=\frac{1}{2\ell_1 \ell_2}
			\left(
			\ell_1+\ell_2+1-\sqrt{\ell_1^2+(2-2\ell_2)\ell_2+\ell_2^2+2\ell_2+1}
			\right).
		\end{equation*}
		There are many choices of critical exponents for which
		$\lambda(\ell_{1l},\ell_{1d})\neq \lambda(\ell_{2l},\ell_{2d})$.
		Hence, by Proposition~22 in \cite{TNB3}, any conjugation between these two maps
		cannot be a H\"{o}lder homeomorphism.
		
		\item Let $X$ be a vector field on $\mathbb{T}^2$, and let $\Phi_X:t\mapsto\Phi_X^t$
		denote its flow on an open set $U$. Let $\psi:\mathbb{T}^2\longrightarrow \mathbb{T}^2$
		be a diffeomorphism. Observe that
		\begin{equation*}
			\left.\frac{d}{dt}\right|_{t=0}\left(\psi\circ\Phi_X^t\circ\psi^{-1}(y)\right)
			=(\psi_*X)_y
			\qquad\text{for every } y\in \psi(U).
		\end{equation*}
		This means that $\psi\circ\Phi_X^t\circ\psi^{-1}$ is the flow of $\psi_*X$.
		Thus, if $X$ generates the map $f$, then $\psi_*X$ generates the map
		\begin{equation*}
			\psi\circ f\circ\psi^{-1}:\psi(\Delta)\approx \mathbb{S}^1\longrightarrow\mathbb{S}^1,\quad
			x\longmapsto \psi\circ \Phi_X^{t(x)}\circ\psi^{-1}(x).
		\end{equation*}
		This ends the proof of Theorem~\ref{Bitheo2}.
	\end{enumerate}
\end{proof}

Let $(X,Y)$ be a pair of Cherry vector fields on $\mathbb{T}^2$, transverse
outside their singular set $\operatorname{Sing}(X,Y)$, and assume that it
generates a map $f_{X,Y}\in\mathscr{L}_c$. Define
$$
\mathcal{C}_{f_{X,Y}}=
\left\{
(X',Y')\in \mathfrak{X}_c^2(\mathbb{T}^2):
f_{X',Y'}=f_{X,Y}
\right\}.
$$

On the punctured torus $\mathring{\mathbb{T}}^2=\mathbb{T}^2\setminus \operatorname{Sing}(X,Y)$,
the vector fields $X$ and $Y$ are nonsingular and transverse. Since
$\mathring{\mathbb{T}}^2$ is an orientable two-dimensional manifold, it admits a symplectic form $\omega$.
Thus, $X$ and $Y$ define two transverse one-dimensional foliations
$\mathcal{F}_X$ and $\mathcal{F}_Y$ on $\mathring{\mathbb{T}}^2$. Since every one-dimensional foliation on a
two-dimensional symplectic manifold is Lagrangian, the triple 
$(\omega,\mathcal{F}_X,\mathcal{F}_Y)$
is a bi-Lagrangian structure on $\mathring{\mathbb{T}}^2$. We denote its Hess
connection by ${}^{X,Y}\nabla$.

Assume that the restriction of this connection to the chosen section
$\mathbb{S}^1$ depends only on $f_{X,Y}$, namely,
$$
\left\{
{}^{X',Y'}\nabla_{\vert \mathbb{S}^1}:
(X',Y')\in \mathcal{C}_{f_{X,Y}}
\right\}
=
\left\{
{}^{X,Y}\nabla_{\vert \mathbb{S}^1}
\right\}.
$$
Then the connection associated with the Cherry map $f_{X,Y}$ is defined by
$$
\nabla^{f_{X,Y}}:={}^{X,Y}\nabla_{\vert \mathbb{S}^1}.
$$

One challenge is to solve the uniqueness problem for $\nabla$. Such a problem
has been addressed in certain nonlinear dynamical systems, see, for instance,
\cite{FAT}.

By combining point 2 of Theorem~\ref{Bitheo2}, Theorem~\ref{c10}
and Lemma~\ref{Bilem1}, we have the following result.
\begin{proposition}
	Let $X, Y \in \mathfrak{X}_c(\mathbb{T})$. Suppose that, for every pair
	$X', Y' \in \mathfrak{X}_c(\mathbb{T})$ satisfying
	$f_{X',Y'} = f_{X,Y}$, there exists a diffeomorphism
	$\psi \in \operatorname{Diff}(\mathbb{T}^2)$ such that $\psi_* X = X'$ and
	$\psi_* Y = Y'$. 
	It follows that
	$$
	{}^{X',Y'}\nabla = {}^{X,Y}\nabla,
	\qquad \forall (X',Y')\in \mathcal{C}_{f_{X,Y}}.
	$$
	Thus, the connection $\nabla^{f_{X,Y}}$ associated with $f_{X,Y}$ is
	well defined; that is, it does not depend on the choice of representative
	$(X,Y)\in \mathcal{C}_{f_{X,Y}}$.
\end{proposition}


\begin{example}\label{Biex1}
	For simplicity of notation, we work on the normalized torus
	$\mathbb{T}^2=\mathbb{R}^2/(2\pi\mathbb{Z})^2$. Consider the vector field
	$$
	X_c(x,y)=\left(c+\cos x+\cos y,\;\sin x+\sin y\right),
	$$ and define $Y_c=R_{\pi/2}X_c\in \mathfrak{X}(\mathbb{T}^2)$, where
	$0<c<2$ and $R_{\pi/2}$ denotes the rotation by angle $\pi/2$ in the tangent
	spaces. Let $a=\arccos\left(-\frac c2\right)$. Since $0<c<2$, we have
	$a\in\left(\frac{\pi}{2},\pi\right)$, and therefore $
	\operatorname{Sing}(X_c)
	=
	\{p_{\operatorname{saddle}}=(-a,a),\; p_{\operatorname{sink}}=(a,-a)\}.
	$
	
	On the section $\Sigma=\{x=0\}$, one has
	$\dot x=c+1+\cos y>0$.
	Hence $\Sigma$ is transverse to the flow of $X_c$.
	Moreover, the rotation number of the first return map associated with $X_c$
	depends on the parameter $c$.
	Thus, by choosing $c$ so that the flow $\Phi_{X_c}^t$ has no periodic orbit,
	the vector field $X_c$ is a Cherry vector field.
	Since $Y_c=R_{\pi/2}X_c$, it follows that $Y_c$ is also a Cherry vector
	field.
	In particular, one may choose $c$ so that the rotation number of
	$f_{X_c,Y_c}$ is the golden mean.	
	Furthermore, $
	\det\left(X_c,Y_c\right)
	=
	\left(c+\cos x+\cos y\right)^2
	+
	\left(\sin x+\sin y\right)^2.
	$
	Consequently, $\det\left(X_c(x,y),Y_c(x,y)\right)=0$ if and only if $X_c(x,y)=0$.
	Therefore,
	$(\omega,\mathcal{F}_{X_c},\mathcal{F}_{Y_c})$ defines an affine bi-Lagrangian
	structure on $\mathring{\mathbb{T}}^2$, where $\omega$ is the canonical
	symplectic form.
\end{example}

\section{Conclusion}
This paper explored a twofold interaction between bi-Lagrangian geometry and
one-dimensional dynamics. We considered a bi-Lagrangian structure
$(\omega,\mathcal{F}_1,\mathcal{F}_2)$ on a manifold $M$ and constructed, in
different ways, lifted bi-Lagrangian structures on the tangent bundle $TM$
and on the cotangent bundle $T^*M$. One construction provided an affine
bi-Lagrangian structure, and, in a particular case, we lifted only the affine
structure; the lifted structure remained affine. We also lifted the natural
action of the symplectomorphism group on the set of bi-Lagrangian structures.

We then applied these ideas to a dynamical setting on the two-dimensional
torus $\mathbb{T}^2$. Since, on $\mathbb{T}^2$, a bi-Lagrangian structure can be
encoded by a pair of transverse nowhere-vanishing vector fields, we used the
canonical action $(\psi,X)\mapsto \psi_*X$ to extend the dynamics to pairs of
Cherry vector fields. Such a pair is induced to yield a bi-Lagrangian structure on the
punctured torus $\mathring{\mathbb{T}}^2$. This extension is then shown to
induce a conjugation action on a subset of Cherry maps, interpreted as circle
maps with a flat interval. Moreover, for certain Cherry maps, linear
connections are defined in terms of their generating bi-Lagrangian data.
It is shown that the maps generated by the same bi-Lagrangian structure do not
belong to the same H\"older homeomorphism class.

The purpose of this construction is to study circle maps with a flat piece via
the linear connections induced by their generating Cherry vector fields. This
may help complete the geometric description of maps in $\mathscr{L}_c$.
Related approaches appear in \cite{GJST,LP1,TNB3}, although the known results
remain incomplete and depend strongly on the critical exponents. For instance,
the case $(\ell_1,\ell_2)=(1.1,2.1)$ is still not fully understood. Finally, it
would be interesting to study the existence of a surjective map
$\mathfrak{X}_c^2(\mathbb{T}^2)\longrightarrow \mathscr{L}_c$,
$(X,Y)\longmapsto f_{X,Y}$.



\small

\begin{thebibliography}{99}
	
	\bibitem{GB1} M.\,N. Boyom.
	The Cohomology of Koszul-Vinberg.
	{\it Pacific Journal of Mathematics} \textbf{225}\,(1) (2006) 119--153. DOI:10.2140/pjm.2006.225.119.
	
	\bibitem{TC} T.\,M. Cherry.
	Analytic Quasi-Periodic Curves of Discontinuous Type on a Torus.
	{\it Proceedings of the London Mathematical Society} \textbf{44}\,(1) (1938) 175--215. DOI:10.1112/plms/s2-44.1.175.
	
	\bibitem{dasilva} A.\,C. da Silva.
	{\it Lectures on Symplectic Geometry}.
	Springer-Verlag, New York, 2006.
	
	\bibitem{1} F. Etayo, R. Santamaria and U. Tr\'{\i}as.
	The geometry of a bi-Lagrangian manifold.
	{\it Differential Geometry and Applications} \textbf{24}\,(1) (2006) 33--59. DOI: 10.1016/j.difgeo.2005.05.004.
	
	\bibitem{FAT} A. Ferrag, A. Aoufi and M. Tadjine.
	Existence and Uniqueness of Solution for Stochastic Nonlocal Random Functional Integral Equation.
	{\it Nonlinear Dynamics and Systems Theory} \textbf{25}\,(1) (2025) 53--58.
	
	\bibitem{GJST} J. Graczyk, L.\,B. Jonker, G. \'{S}wi\c{c}tek, F.\,M. Tangerman and J.\,J.\,P. Veerman.
	Differentiable Circle Maps with a Flat Interval.
	{\it Communications in Mathematical Physics} \textbf{173}\,(3) (1995) 599--622. DOI: 10.1007/BF02101658
	
	\bibitem{7} H. Hess.
	Connections on symplectic manifolds and geometric quantization.
	In: {\it Differential Geometric Methods in Mathematical Physics}.
	(P.\,L. Garc{\'{\i}}a, A. P{\'e}rez-Rend{\'o}n and J.\,M. Souriau).
	Springer-Verlag, Berlin, 1980, 153--166. DOI: 10.1007/BFb0089731
	
%
	\bibitem{LP1} L. Palmisano.
	 A Phase Transition for Circle Maps and Cherry Flows.
	{\it Communications in Mathematical Physics} \textbf{321}\,(1) (2013) 135--155. DOI:  10.1007/s00220-013-1685-2.
	
	\bibitem{TNB2} B. Tangue Ndawa.
	Infinite Lifting of an Action of Symplectomorphism Group on the set of bi-Lagrangian structures.
	{\it Journal of Geometry and Mechanics} \textbf{14}\,(3) (2022) 409--426. DOI: 10.3934/jgm.2022006.
	
	\bibitem{TNB3} B. Tangue Ndawa.
	Rigidity of Fibonacci Circle Maps with a Flat Piece and Different Critical Exponents.
	{\it Journal of Dynamics and Control Systems} \textbf{31}\,(13) (2025) 1--42. DOI: 10.1007/s10883-024-09687-z
	
	\bibitem{YKI} K. Yano and S. Kobayashi.
	Prolongations of tensor fields and connections to tangent bundles I. 
	{\it Journal of the Mathematical Society of Japan} \textbf{18}\,(2) (1966) 195--210. DOI: 10.2969/jmsj/01820194.
	
	\bibitem{YKII} K. Yano and S. Kobayashi.
	Prolongations of tensor fields and connections to tangent bundles II. 
	{\it Journal of the Mathematical Society of Japan} \textbf{18}\,(3) (1966) 236--246. DOI: 10.2969/jmsj/01830236.
	
\end{thebibliography}
\end{document}